\documentclass{article}

\usepackage{amsthm}
\usepackage{amsfonts}
\usepackage{amssymb}
\usepackage{amsmath}
\usepackage{tikz}
\usepackage{acronym}
\usepackage{hyperref}
\hypersetup{bookmarksdepth=section}
\usepackage{cleveref}
\usepackage{subcaption}
\usepackage{graphicx}
\usepackage{url}
\usetikzlibrary{positioning, shapes.geometric, calc}

\DeclareMathOperator{\lv}{lev}
\DeclareMathOperator{\lc}{lc}
\DeclareMathOperator{\win}{win}
\newcommand{\keywords}[1]{\par \noindent \textbf{Keywords:} #1}

\begin{document}
\title{The Construction of Directed Partial Orders on the Complex Field by Puiseux Series and Stratification}
\author{Zitian Zhu \thanks{Corresponding author}}
\date{\today}

\theoremstyle{definition}
\newtheorem{theorem}{Theorem}[section]
\newtheorem{lemma}[theorem]{Lemma}
\newtheorem{proposition}[theorem]{Proposition}
\newtheorem{property}[theorem]{Property}
\newtheorem{corollary}[theorem]{Corollary}
\newtheorem{definition}{Definition}[section]
\newtheorem{remark}{Remark}[section]
\newtheorem{example}{Example}
\newtheorem{exercise}{Exercise}

\maketitle
\begin{abstract}
\par Birkhoff and Pierce asked in 1956 whether the complex field $\mathbb{C}$ admits a lattice order compatible with the ring operations. Schwartz and Yang proved in 2011 that $\mathbb{C}$ admits directed partial orders. We give a new construction that by using a transcendence basis of $\mathbb{C}$ over $\mathbb{Q}$, we embed $\mathbb{C}$ into a field of Puiseux series $F\{\{s\}\}$, so that every complex number acquires a series expansion, and we define positivity by a condition on the leading coefficient. Every directed partial order we construct is placed in a stratification by how much of the expansion a rule has to read. Reading one term, a window of fixed width, or a sensitivity condition on small changes, we obtain five nested conditions $A\Rightarrow B\Rightarrow C\Rightarrow D\Rightarrow E$ and six types. Relative to the fixed level valuation, cones satisfying $B$ correspond to the cards satisfying a product and a sum condition. Cones satisfying $C$ are never lattice orders, which excludes Types $1,2,3$. We give an example of Type $4$ that is not a lattice order either. The existence of a lattice order on $\mathbb{C}$ and the existence of orders of Types $5$ and $6$ still remain open. This is a stratification rather than a complete classification because only Types $1$--$4$ are shown to be nonempty.
\end{abstract}
\keywords{directed partial order, cone, Puiseux series, valuation, lattice-ordered ring, Birkhoff--Pierce problem}
\newpage
\tableofcontents

\section{Introduction}
\par In 1956, Birkhoff and Pierce asked the question whether the complex field $\mathbb{C}$ admits a lattice order compatible under the rings \cite{bp56}; in 2011, Schwartz and Yang proved that the directed partial order are admitted on $\mathbb{C}$ \cite{sy11} and later found that a field admits one if and only if it has characteristic $0$ and is real or has the transcendence degree at least 1 over $\mathbb{Q}$ \cite{sy23}. Furthermore, for a non-archimedean totally ordered field $F$, the directed partial orders on $F(i)$ compatible with the order of $F$ were studied and classified by Rump and Yang \cite{ry14}, Ma, Wu, and Zhang \cite{mwz17} \cite{mwz18} \cite{mwz19}, and Xu and Zhang \cite{xz20}. Recently, Wang, Yuan, Zhang, and Zhu \cite{wyzz26} also give a related example by construction, stating the conjecture that there does not exist lattice order on the complex field. None of these orders is a lattice order. However, He \cite{z26} gives a certain answer that $\mathbb{C}$ admits at least $2^{\aleph_{0}}$ distinct compatible lattice orders.
\par This paper has two parts. The first part is a direct construction on $\mathbb{C}$ as an abstract field without ordered subfield given in advance. We write $\mathbb{C}$ as the algebraic closure of $F(t)$, where $t$ is transcendental and $F$ is an algebraically closed subfield, put $s=1/t$, and embed $\mathbb{C}$ into the field of Puiseux series over $F$ in $s$. Thus, every complex number becomes a series like a decimal expansion. The second part constructs a stratification to sort all directed partial orders on $\mathbb{C}$. We sort the orders by how much of the series must be read. With the restriction, we introduce five nested conditions and show every cone falls into one of the six classes exactly. The first four classes are proved and described to be nonempty, but the last two still remain some possibilities for their existence. The conditions are shaped by relative infinitesimals, the pairs of ring and ideal, and the leading coefficients. However, some of the types have no example and are unknown for the existence of the lattice order. The Birkhoff-Pierce problem still remains open in this paper.
\par \textit{Section} \ref{pre} includes the necessary concepts and terms that ones need to know for the paper. \textit{Section} \ref{con} constructs the Puiseux expansion of complex numbers and the basic orders, proving the basic orders. \textit{Section} \ref{st} defines the six types of the stratification of the directed partial orders on $\mathbb{C}$, demonstrating the Covering Lemma and explaining each type of directed partial order.
\section{Preliminaries} \label{pre}
\par The purpose of this section is to mention some concepts and notations that appear in the proofs. For a subset $I$ of $\mathbb{C}$, we write $I^{\times} = I \backslash \lbrace 0 \rbrace$. We define $i = \sqrt{-1}$.
\par Let $K\subseteq L$ be fields. We use two standard facts: the set of elements of an algebraically closed field that are algebraic over a given subfield is an algebraically closed subfield; if $L$ is algebraic over $K$, then every embedding of $K$ into an algebraically closed field extends to an embedding of $L$. When we decide the variable and the coefficients, three things guarantees the choice. First, the variable $t$ must be transcendental over the field $F$ of coefficients. Second, $F$ must be algebraically closed because of the Newton-Puiseux Theorem \ref{thm:NP}. Third, $F$ must be large enough in order that every complex number is algebraic over $F(t)$ in order to let $\mathbb{C}$ inside the series field. A transcendence basis $T$ is the tool to achieve them at the same time.  We fix a transcendence basis $T$ of $\mathbb{C}$ over $\mathbb{Q}$ and one element $t \in T$. Put $s = \frac{1}{t}$, and let $F$ be the set of complex numbers which are algebraic over the field $\mathbb{Q}(T \backslash \lbrace t \rbrace)$. The elements of $F$ are the constants.
\begin{theorem}[Newton-Puiseux]\label{thm:NP}
If $k$ is an algebraically closed field of characteristic $0$, then $k \lbrace \lbrace s \rbrace \rbrace$ is
algebraically closed.
\end{theorem}
\par The theorem is used for Lemma \ref{lem32} and involved many unnecessary concepts for this paper. We just use it as a standard conclusion. Next, we introduce the crucial terms.
\begin{definition}
A \textit{cone} on a field $K$ is a subset $P$ of $K$ with
\[
0 \in P, \qquad P+P\subseteq P,\qquad P\cdot P\subseteq P,\qquad P\cap(-P)=\{0\}.
\]
The elements of $P^{\times}$ are called \textit{positive}.
\end{definition}
 
A partial order on a field and a cone are the same thing. If $\le$ is a partial order on the field $K$, then the set $P=\{x\in K:x\ge0\}$ is a cone. Conversely, if $P$ is a cone, then $x\le y\quad\Longleftrightarrow\quad y-x\in P$ is a partial order on the field $K$ whose set of elements $\ge0$ is $P$. A partial order is \textit{directed} if for every $x,y$ there is $z$ with $x\le z$ and $y\le z$; equivalently, $P-P=\mathbb{C}$. Therefore, we say that a cone $P$ is a total order, a lattice order or a directed partial order when the order it defines is one. We shall mostly take the field as $\mathbb{C}$ or sometimes a subfield of $\mathbb{C}$. 
\par These two lemmas are used in the later part:
\begin{lemma}\label{lem21}
Let $P$ be a cone on $\mathbb{C}$ and $p\in P^{\times}$. Then neither $ip$ nor $-ip$ lies in $P$.
\end{lemma}
\begin{proof}
If $ip\in P$, then $(ip)^2=-p^2\in P$. Also, $p^2\in P$. Thus, $p^2\in P\cap(-P)=\{0\}$, which is false. For $-ip$, it is the same thing.
\end{proof}
 
\begin{lemma}[Wiggle lemma]\label{lemwiggle}
Let $P$ be a cone on $\mathbb{C}$ and $x,y\in\mathbb{C}$ be incomparable. Suppose $h=x\vee y$ exists. Put $p=h-x$ and $q=h-y$. Then $p,q\in P^{\times}$ and $[-p,p]\cap[-q,q]=\lbrace 0 \rbrace$.
\end{lemma}
\begin{proof}
Because $h$ is an upper bound, $p,q\in P$. The value $p=0$ forces the relation $ y \leq x$. Let $d\in[-p,p]\cap[-q,q]$. Then, $(h-d)-x=p-d\in P$ and $(h-d)-y=q-d\in P$, so $h-d$ is an upper bound of $x$ and $y$. Since $h$ is the least upper bound, $h\le h-d$, $-d\in P$. Similarly, $h+d$ is an upper bound, so $d\in P$. Thus, $d=0$.
\end{proof}
 
\begin{definition}
A \textit{valuation} on $\mathbb{C}$ is a surjective map $w\colon\mathbb{C}^{\times}\to\Gamma_w$ onto a totally ordered abelian group such that $w(xy)=w(x)+w(y)$, $w(x+y)\ge\min\{w(x),w(y)\}$.
It is \textit{nontrivial} if $\Gamma_w \ne0$. We put $\mathcal{O}_w=\{x:w(x)\ge0\}\cup\{0\}$, $\mathfrak{m}_w=\{x:w(x)>0\}\cup\{0\}$, $I^w_\delta=\{x:w(x)>\delta\}\cup\{0\}$ for $\delta\in\Gamma_w$ with $\delta\ge 0$.
\end{definition}

\par To be noticed, when the valuation is the level valuation $\lv$ of Section \ref{con}, we write $\mathcal{O} = \mathcal{O}_{\lv}$, $\mathfrak{m} = \mathfrak{m}_{\lv}$, and $I_{\delta} = I_{\delta}^{\lv}$ for $\delta \geq 0$. For a general nontrivial valuation $w$, we keep the notation $I_{\delta}^{w}$.

\begin{lemma}\label{lem:valring}
Let $w$ be a nontrivial valuation on $\mathbb{C}$.
\begin{enumerate}
\item $\mathcal{O}_w$ is a subring of $\mathbb{C}$ with quotient field $\mathbb{C}$, and for every $x\in\mathbb{C}^{\times}$ at least one of $x$, $x^{-1}$ lies in $\mathcal{O}_w$.
\item Each $I^w_\delta$ is a nonzero ideal of $\mathcal{O}_w$, $I^w_0=\mathfrak{m}$, and $1+I^w_\delta$ is a group under multiplication.
\end{enumerate}
\end{lemma}
\begin{proof}
(a) $\mathcal{O}_w$ is closed under sums and products by the two rules. Since $w(x^{-1})=-w(x)$, one of $x$ and $x^{-1}$ has the value at least $0$. Thus, every $x$ is a quotient of
elements of $\mathcal{O}_w$. (b)  $I^w_\delta\ne0$ because $\Gamma_w\ne0$ has elements larger than $\delta$. For $\varepsilon,\varepsilon'\in I^w_\delta$, we have $(1+\varepsilon)(1+\varepsilon')\in1+I^w_\delta$ and $(1+\varepsilon)^{-1}=1 - \varepsilon(1+\varepsilon)^{-1}$ with $w(\varepsilon(1+\varepsilon)^{-1})=w(\varepsilon)$.
\end{proof}
\par Let $k$ be a field. A \textit{Laurent series} over $k$ in a variable $s$ is a formal sum $\sum_{m\ge m_0}c_{m}s^m$ with $m_{0}\in\mathbb{Z}$ and $c_m\in k$. Laurent series are added and multiplied like polynomials and form a field $k((s))$ which contains the field $k(s)$ of rational functions. For an integer $n\ge1$, let $k((s^{\frac{1}{n}}))$ be the field of Laurent series in the variable $s^{\frac{1}{n}}$. If $n$ divides $m$, then $k((s^{\frac{1}{n}}))\subseteq k((s^{\frac{1}{m}}))$.

\section{Complex Numbers with Puiseux Series} \label{con}
\par In this section, we rewrite every complex number as the form of series and construct the basic directed partial orders on $\mathbb{C}$. There are three steps: we isolate one transcendental number $t$, embed $\mathbb{C}$ into a field of Puiseux series in $s=\frac{1}{t}$, and define positivity by a condition on the first coefficient.
\begin{definition}\label{def:puiseux}
The field of \textit{Puiseux series} over $k$ in $s$ is
\[
k \lbrace\lbrace s \rbrace \rbrace = \bigcup_{n\ge1}k((s^{\frac{1}{n}})).
\]
\end{definition}
\par Thus, a nonzero Puiseux series is a formal sum $\sum_qc_qs^q$ with $c_q\in k$, whose exponents $q$ are rational numbers with a common denominator and are bounded below. Fractional exponents cannot be avoided. A square root of $s$ is not a Laurent series in $s$. 
\begin{lemma} \label{lem32}
    There is an embedding of fields $\iota\colon\mathbb{C}\to F\lbrace \lbrace s \rbrace \rbrace$ which is the identity on $F$ and sends $t$ to $s^{-1}$.
\end{lemma}
\begin{proof}
    The set of elements of $\mathbb{C}$ that are algebraic over a subfield is an algebraically closed subfield, and it contains every number algebraic over $\mathbb{Q}$. If $t$ were algebraic over $F$, it would be algebraic over $\mathbb{Q}(T\setminus\lbrace t\rbrace)$, contradicting the algebraic independence of $T$. Thus, $F$ which contains $\Bar{\mathbb{Q}}$ is an algebraically closed subfield of $\mathbb{C}$, and $t$ is transcendental over $F$. Therefore, there is an isomorphism of $F(t)$ onto the subfield $F(s)$ of $F \lbrace \lbrace s \rbrace\rbrace$ which is the identity on $F$ and sends $t$ to $s^{-1}$. By Theorem \ref{thm:NP}, the field $F \lbrace \lbrace s \rbrace\rbrace$ is algebraically closed. Because every complex number is algebraic over $\mathbb{Q}(T)$ and $\mathbb{Q}(T) \subseteq F(t)$, every complex number is algebraic over $F(t) = F(s)$. Hence, $\mathbb{C}$ is algebraic over $F(t)$. Thus, the isomorphism generalizes to an embedding of $\mathbb{C}$.
\end{proof}
\subsection{The Expansion}
\par Now, $\mathbb{C}$ and the field of Puiseux series are two different objects. An order defined by "the first coefficient" needs every complex number to be a series. The embedding of Lemma \ref{lem32} provides that it is a dictionary to translate each complex number into a series without changing sums and products. After the translation, the cone conditions become the statements about adding and multiplying series.
\par We fix an embedding $\iota$ as the one in Lemma \ref{lem32} and identify a complex number $z$ with the series $\iota(z)$. Thus, every nonzero $z \in \mathbb{C}$ has an expansion \[
    z=\sum_{q}z_q\,s^q, \; z_q\in F .
\]
\par It is like a decimal expansion. The image of $\mathbb{C}$ is an algebraically closed field containing $F$ and $s$, so, for every $n \geq 1$, it contains all $n$ roots of the polynomial $x^{n}-s$ existing in $F \lbrace \lbrace s \rbrace \rbrace$. The series $s^{\frac{1}{n}}$ is one of them. Thus, $s^{q}$ and $t^{q}$ are complex numbers for every rational number $q$.
\begin{example} \label{ex1}
\par Let $z$ be a complex number with $z^2=t^2+t$. Only algebraic over $F(t)$, it is not a rational function of $t$. Writing $t^2+t=t^2(1+s)$ and expanding the square root by the binomial series gives
\[
z=s^{-1}(1+s)^{1/2}=s^{-1}+\tfrac12-\tfrac18\,s+\cdots.
\]
Thus, $z$ has the size of $s^{-1} = t$ and the first coefficient $1$. An order reading the first term only calls $z$ positive. The other solution $-z$ has the series $-s^{-1}-\tfrac12+\cdots$, with first coefficient $-1$. Decided by the embedding $\iota$, $t$ determines which solution begins with $+s^{-1}$.
\par Another simple example is a square root of $t$, which is $s^{-\frac{1}{2}}$, explaining why fraction exponents are needed.
\end{example}
\par The expansion is the way of writing the number in order that the orders below can read its coefficients. By the axiom of choice, it depends on the choice of $T$ and of the embedding $\iota$. We know no directed partial order on $\mathbb{C}$ is defined without this choice.

\subsection{Level, Leading Coefficient, and Normalized Series}
\par A rule reading the beginning of a series observes the things of the place where the series begins and with which coefficient. The first one is defined as the level which measures the size of the number and replaces the absolute value. It is useless here because it does not come from an order on $\mathbb{C}$. The second one is the leading coefficient. 
\begin{definition}\label{def:level}
Let $z \in \mathbb{C}$ be nonzero.
\begin{enumerate}
\item The \textit{level} $\lv(z)$ is the smallest exponent $q$ with $z_q\ne0$.
\item The \textit{leading coefficient} is $\lc(z)=z_{\lv(z)}\in F^{\times}$.
\item The \textit{normalized series} is $u=z/s^{\lv(z)}=\sum_{r\ge0}u_rs^r$, so
$u_0=\lc(z)$.
\end{enumerate}
\end{definition}
\par The property is obtained directly from the definition: $\lv(c) = 0$, $\lc(c)=c$ $c \in F^{\times}$; $\lv(s^{q}) = q$, $\lc(s^{q})=1$ $q \in \mathbb{Q}$; $\lc (-z) = - \lc (z)$. The normalized series removes the size, so the numbers of different sizes can be compared digit by digit. To be noticed, a very negative level implies a very large number due to $s^{-1} = t$. For instance, in Example \ref{ex1}, the level is $-1$; the leading coefficient is $1$; the normalized series is $(1+s)^{\frac{1}{2}} = 1 + \frac{1}{2}s - \frac{1}{8}s^{2} + \cdots$.
\begin{lemma} \label{lem33}
    The level is a valuation on $\mathbb{C}$ with value group $\mathbb{Q}$.
\[
\lv(zw)=\lv(z)+\lv(w),\quad\lv(z+w)\ge\min\{\lv(z),\lv(w)\}.
\]
For nonzero $z,w$, we have the following.
\begin{enumerate}
\item $\lc(zw)=\lc(z)\lc(w)$;
\item $\lv(z)<\lv(w)$ $\Rightarrow$ ($\lv(z+w)=\lv(z)$, $\lc(z+w)=\lc(z)$);
\item $\lv(z)=\lv(w)$,  $\lc(z)+\lc(w)\ne0$ $\Rightarrow$ ($\lv(z+w)=\lv(z)$, $\lc(z+w)=\lc(z)+\lc(w)$).
\end{enumerate}
\end{lemma}
\begin{proof}
\par Rewrite $z=s^q(a+\cdots)$ and $w=s^{q'}(b+\cdots)$ with $a=\lc(z)$, $b=\lc(w)$ where the dots are terms with positive exponents. Then, $zw=s^{q+q'}(ab+\cdots)$, and $ab\ne0$. If $q<q'$, then $z+w=s^q(a+\cdots)$. If $q=q'$, then $z+w=s^q((a+b)+\cdots)$, which has level $q$ and leading coefficient $a+b$ when $a+b\ne0$. Otherwise, it appears a higher level. The value group is $\mathbb{Q}$ because $s^q\in\mathbb{C}$, $\forall q$.
\end{proof}
\par Lemma \ref{lem33} explains how the two data are the right ones. Under multiplication and addition, they change in a predictable way, and this is all to verify the cone conditions. Furthermore, we write $\mathcal{O}$ for the ring of numbers of level $\ge0$ (together with $0$), $\mathfrak{m}$ for the ideal of numbers of level greater than $0$, and $I_\delta$ for the ideal of numbers of level greater than $\delta$.
\subsection{The Basic Order by the First Coefficient}
\par We call $z$ positive when its leading coefficient is positive. However, the leading coefficient lies in $F$ which contains $i$ and has no total order, so the word "positive" has no actual meaning yet. The necessary thing is a set $S$ of admissible leading coefficients, and the conditions on $S$ are dictated by the cone conditions for $\mathbb{C}$. The set $S$ must be closed under multiplication and addition. The first terms cannot be canceled by the condition $S\cap(-S)=\{0\}$. These are exactly the conditions for a cone on $F$. Thus, a cone on the small field $F$ becomes to a cone on $\mathbb{C}$.
\begin{theorem} \label{thm:basic}
Let $S$ be a cone on the field $F$ of constants with $1\in S$, and $d(S)=\{0\}\cup\{z\in\mathbb{C}^{\times}:\ \lc(z)\in S\}.$ Then $d(S)$ is a directed partial order on $\mathbb{C}$ containing $1$, and $d(S)\cap F=S$.
\end{theorem}
\begin{proof}
\par By $\lc(1)=1$, we have $1\in d(S)$. A nonzero constant $c$ has $\lc(c)=c$, so $d(S)\cap F=S$. Let $z,w\in d(S)$ be nonzero. By Lemma \ref{lem33}, for multiplication, $\lc(zw)=\lc(z)\lc(w)\in S\cdot S\subseteq S$; for addition, if $\lv(z)<\lv(w)$, then $\lc(z+w)=\lc(z)\in S$ and the case $\lv(z)>\lv(w)$ is the same. If $\lv(z)=\lv(w)$, then $\lc(z)+\lc(w)\in S$ and it is not $0$, or $\lc(z)=-\lc(w)\in S\cap(-S)=\{0\}$. Thus, $\lc(z+w)=\lc(z)+\lc(w)\in S$. The condition $S\cap(-S)=\{0\}$ is to prevent the first terms to be canceled.
\par No element is positive and negative at the same time. If $z\ne0$ and $z,-z\in d(S)$, then $\lc(z)$ and $-\lc(z)$ lie in $S$, which is impossible.
\par Let $z\in\mathbb{C}$ be nonzero. We choose an integer $N>-\lv(z)$ and let $w=t^N=s^{-N}$. Because $\lc(w)=1\in S$ and $\lv(w)=-N<\lv(z)$, $\lc(z+w)=\lc(w)=1$, and $w$ and $z+w$ lie in $d(S)$ with $z=(z+w)-w$, also $0=w-w$. Thus, the order is directed. Lastly, by Lemma \ref{lem21}, no cone on $\mathbb{C}$ is total, so the order is not total.
\end{proof}
\par We define For $D = d(\mathbb{Q}_{\geq 0})=\{0\}\cup\{z\in\mathbb{C}^{\times}:\ \lc(z)\in \mathbb{Q}_{\geq 0}\}$. For $D$, the number $z=\frac{3}{2}t^2-it+\sqrt2=\frac{3}{2}s^{-2}-is^{-1}+\sqrt2$ is positive, because its first coefficient $\frac{3}{2}$ is a positive rational number. The later coefficients $-i$ and $\sqrt2$ are not looked at. The number $\sqrt2\,t$ and $-\sqrt2\,t$ is not positive because the first coefficient $\pm\sqrt2$ is not rational. Every number $z$ is below $t^N$ for $N$ large, and this is how the powers of $t$ make the order directed. The following corollary is the case of the rational number field.
\begin{corollary}
$D\cap F=\mathbb{Q}_{\ge0}$. In particular, $D\cap\Bar{\mathbb{Q}}=\mathbb{Q}_{\ge0}$.
\end{corollary}
\begin{proof}
\par By Theorem \ref{thm:basic}, we substitute $S=\mathbb{Q}_{\ge0}$ and can have that $D=d(\mathbb{Q}_{\ge0})=\{0\}\cup\{z\in\mathbb{C}^{\times}:\ \lc(z)\in\mathbb{Q}_{>0}\}$ is a directed partial order on $\mathbb{C}$, and it is not total. Thus, $D \cap F = \mathbb{Q}_{\geq 0}$. In addition, by the proof of Lemma \ref{lem32}, $\Bar{\mathbb{Q}} \subseteq F$.
\end{proof}
\par The construction involves four choices of the set $S$, the embedding $\iota$, the element $t$, and the basis $T$. The next example changes them and shows that each change produces a different order. This also contributes to the stratification in Section \ref{st} because there is no single Puiseux order on $\mathbb{C}$ and the description of all directed partial orders must adapt for this freedom.
\begin{example} \label{ex:many}
\begin{itemize}
\item We can make the sets of first coefficients smaller: $d(\mathbb{N})\subsetneq d(\mathbb{Q}_{\ge0})=D$.
\item Let $L$ be a subfield of $F$ with a total order, and $L_{\ge0}$ its cone. Then $d(L_{\ge0})$ is a directed partial order with $d(L_{\ge0})\cap F=L_{\ge0}$. There are many such $L$. For a subset $A$ of $T\setminus\{t\}$, the field $\mathbb{Q}(A)$ embeds in $\mathbb{R}$ by sending $A$ to algebraically independent real numbers. Let $S_A$ be the set of elements of $\mathbb{Q}(A)$ whose image is \(\ge0\). Then, $S_A-S_A=\mathbb{Q}(A)$, and $\mathbb{Q}(A)\ne\mathbb{Q}(B)$ for $A\ne B$ because $T$ is algebraically independent. Thus, the orders $d(S_A)$ are pairwise distinct, and there are at least as many of them as there are subsets of $T\setminus\{t\}$.
\item For $q\in\mathbb{Q}$, let $\chi(q)=e^{2\pi iq}$, a root of unity and an element of $F$. The map $\sigma\bigl(\sum_qc_qs^q\bigr)=\sum_q\chi(q)c_qs^q$ is an automorphism of $F\{\{s\}\}$ because $\chi(q+q')=\chi(q)\chi(q')$. It is the identity on $F((s))$ because $\chi(m)=1$ for $m\in\mathbb{Z}$. Thus, $\sigma\circ\iota$ is another embedding as in Lemma \ref{lem32}. Its leading coefficient is $\chi(\lv(z))\lc(z)$, and its basic order is $D^\sigma=\{0\}\cup\{z\in\mathbb{C}^{\times}:\ \chi(\lv(z))\lc(z)\in\mathbb{Q}_{>0}\}$. The number $s^{\frac{1}{2}}$ lies in $D$ but not in $D^\sigma$ since $\chi(\frac{1}{2})=-1$. Therefore, the
order depends on the embedding.
\item Let $t'\in T$ with $t'\ne t$. Let $D'$ be the order obtained in the same way from $t'$. Then, $t'\in D'$. However, $t'$ is a constant for the expansion in $s$ and $\lc(t')=t'\notin\mathbb{Q}$, so $t'\notin D$. Thus, $D'\ne D$.
\end{itemize}
\end{example}
\par The Birkhoff-Pierce problem asks for a lattice order, so the first question about the new order is whether it is the one. The natural pair usually for testing is $0$ and $i$, because these two numbers are incomparable on $\mathbb{C}$ that is proved in Lemma \ref{lem21}. The idea of the test is that an upper bound of $0$ and $i$ must be a "large" number (of negative
level). However, among large numbers,  there is no the smallest one because from any large number one can subtract a number that is still large, and the result is an upper bound. The following lemma and theorem demonstrate the idea more clearly.
\begin{lemma} \label{lem:ub}
Let $S$ be a cone on $F$ with $1\in S$ contained in a subfield $L$ of $F$ with
$i\notin L$. Then $z\in\mathbb{C}$ is an upper bound of $0$ and $i$ for the order $d(S)$ if and only if $\lv(z)<0$ and $\lc(z)\in S$.
\end{lemma}
\begin{proof}
\par $z$ is an upper bound of $0$ and $i$ if and only if $z\in d(S)$ and $z-i\in d(S)$. Notice that $\lv(i)=0$, $\lc(i)=i$ and $\pm i\notin L\supseteq S$.
\par Let $z$ be an upper bound. Then, $z\ne0$, since $\lc(-i)=-i\notin S$. Therefore, $\lc(z)\in S$. If $\lv(z)>0$, then $\lc(z-i)=-i\notin S$ by Lemma \ref{lem33}. If $\lv(z)=0$, then $c=\lc(z)\in L$, so $c-i\ne0$ and $\lc(z-i)=c-i\notin L$ by Lemma \ref{lem33}. Both of them never exist, so $\lv(z)<0$. Conversely, if $\lv(z)<0=\lv(i)$ and $\lc(z)\in S$, then $z\in d(S)$ and $\lc(z-i)=\lc(z)\in S$, so $z-i\in d(S)$.
\end{proof}
\begin{theorem} \label{thmnol}
Let $S$ be a cone on $F$ with $1\in S$ which is contained in a subfield of $F$ not containing $i$. Then, $d(S)$ is not a lattice order.
\end{theorem}
\begin{proof}
\par A totally ordered field does not contain $i$, so the examples are covered. Suppose $h=0\vee i$ exists. By Lemma \ref{lem:ub}, $a=-\lv(h)>0$ and $\lc(h)\in S$. Put $h'=h-t^{\frac{a}{2}}$. Since $\lv(t^{\frac{a}{2}})=-\frac{a}{2}>\lv(h)$, Lemma \ref{lem33} gives $\lv(h')=\lv(h)<0$ and $\lc(h')=\lc(h)\in S$. Thus, $h'$ is an upper bound of $0$ and $i$ by Lemma \ref{lem:ub}. However, $h-h'=t^{\frac{a}{2}}$ is positive, so $h'<h$, contradicting the choice of $h$.
\end{proof}
\par This theorem is a motivating special case of Theorem \ref{thm:nolattice}.
\par The proof shows that the set of upper bounds of $0$ and $i$, which contains $t$, has no minimal element at all. Nevertheless, the later section removes the hypothesis on $S$ and extends the statement to all orders of some types of directed partial orders. Section \ref{con} uses only that $\mathbb{C}$ is an algebraically closed field of characteristic $0$ with a transcendence basis $T\ne\emptyset$ over $\mathbb{Q}$. Therefore, Theorems \ref{thm:basic} and \ref{thmnol} hold for every algebraically closed field of characteristic $0$ and transcendence degree at least $1$ over $\mathbb{Q}$. The method cannot work for the field $\Bar{\mathbb{Q}}$ of algebraic numbers, which has no transcendental element to have positivity. Indeed, $\Bar{\mathbb{Q}}$ has no directed partial order at all \cite{sy23}.
\section{The Stratification on the Complex Field} \label{st}
\par Section \ref{con} produced some examples by the expansion of complex numbers with the Puiseux series. In this section, we attempt to stratify all the directed partial order in $\mathbb{C}$, including those not given by any series. Therefore, we need a property valid for an arbitrary cone and separating the orders we know. The proof of Theorem \ref{thmnol} suggests such a property. The basic orders fail to be lattice orders because a positive number still stays positive when it is changed a little, causing an upper bound can always be lowered a little. A lattice order must behave differently when it exists. Hence, we stratify the orders by their sensitivity to small changes.
\subsection{Windows}
\par The basic orders read one term of the series. To find more orders and stratify them, we let a rule read a longer piece by inducing the definition of the window in order to make "the first few terms" clear. The window is dependent on the normalized series, so its width is measured from the first term of each number and not from a fixed exponent, allowing the windows to be compatible with multiplication.  Because the width $\delta$ is fixed, the rule reads the same amount of every number. 
\begin{definition}\label{def:window}
Fix a rational number $\delta\ge0$, the \textit{width}.
\begin{enumerate}
\item A \textit{truncated series} of width $\delta$ is a finite sum
$a=\sum_{0\le r\le\delta}a_rs^r$ with $a_r\in F$ and rational exponents.
\item For $0\le\rho\le\delta$ and $b\in T_{\delta}$, the element $s^{\rho} b\in T_{\delta}$ is $b$ shifted by $\rho$ with the terms beyond $\delta$ deleted.
\item The \textit{window} of a nonzero $z\in\mathbb{C}$ is its normalized series cut off after the exponent $\delta$: $\win(z)=\sum_{0\le r\le\delta}u_rs^r\in T_{\delta}^{\times}$.
\end{enumerate}
\end{definition}
\par Truncated series are added in terms and multiplied by multiplying out and deleting all terms of exponent larger than $\delta$. They form a ring $T_{\delta}$. An element of $T_{\delta}$ is a unit when its constant term $a_0$ is nonzero. $T_{\delta}^{\times}$ denotes the set of units. There are also some examples about the window: For $\delta=0$, the window is the leading coefficient. For $\delta=1$ and the expansion $z=3s^{-2}+5s^{-3/2}-s^{-1}+7s^{-1/2}+\cdots$, we have $\lv(z)=-2$ and $\win(z)=3+5s^{1/2}-s$.
\par The following lemma translates the same level and window into equal up to a factor $1+\varepsilon$ with $\varepsilon$ small, later used in the types for classification.
\begin{lemma}\label{lem:window}
\begin{enumerate}
\item $\win(zw)=\win(z)\win(w)$.
\item Every pair $(q,a)$ with $q\in\mathbb{Q}$ and $a\in T_{\delta}^{\times}$ is the level and window of some complex number.
\item Two nonzero numbers $z,z'$ have the same level and window if and only if $z'=z(1+\varepsilon)$ with $\varepsilon\in I_\delta$.
\item Let $\lv(z)=q$, $\lv(z')=q+\rho$ with $\rho\ge0$, $a=\win(z)$, $b=\win(z')$. If $\rho>\delta$, then $z+z'$ has level $q$ and window $a$. If $\rho\le\delta$ and the constant term of $a+s^\rho b$ is nonzero, then $z+z'$ has level $q$ and window $a+s^\rho b$.
\end{enumerate}
\end{lemma}
\begin{proof}
\begin{enumerate}
    \item The normalized series of $zw$ is the product of the two normalized series, and terms of exponent larger than $\delta$  contribute to exponents larger than $\delta$.
    \item A truncated series is a finite sum of terms $a_rs^r$ with $s^r\in\mathbb{C}$.
    \item If \(z,z'\) have the same level and window, then \(z'/z\) has the level of \(0\), and \(z'/z-1\) has the level larger than \(\delta\); hence \(z'=z(1+\varepsilon)\) with \(\varepsilon\in I_\delta\). Conversely, if \(z'=z(1+\varepsilon)\) with \(\varepsilon\in I_\delta\), then by (1) and \(\win(1+\varepsilon)=1\), we have \(\win(z')=\win(z)\win(1+\varepsilon)=\win(z)\).
    \item Let $z=s^qu$ and $z'=s^{q+\rho}u'$. Then, $z+z'=s^q(u+s^\rho u')$, and $u+s^\rho u'$ agrees with $a+s^\rho b$ up to exponent $\delta$.
\end{enumerate}
\end{proof}

\subsection{Rules and Small Change}
\par Let $P$ be a cone on $\mathbb{C}$. For the basic order $D= d(\mathbb{Q}_{\ge0})$ in Theorem \ref{thm:basic}, a number is positive when its leading coefficient is a positive rational number. To be noticed, multiplying a positive number $z$ by $1+\varepsilon$ with $\lv(\varepsilon)>0$ does not change the leading coefficient, so $z(1+\varepsilon)$ is positive. The classification is built on the observation that an order is described by the small changes $z\mapsto z(1+\varepsilon)$ which can never destruct positivity. The word "small" means that it is specified by a pair.
\par For an arbitrary cone with no given expansion, we need a notion of smallness using only in the field operations. A set $\mathfrak{a}$ of small numbers should have two properties that the sum of two small numbers is small and a small number times a moderate number is small. If $R$ denotes the set of moderate numbers, this illustrates that $\mathfrak{a}$ is an ideal of the ring $R$. We ask that the quotient field of $R$ is $\mathbb{C}$ so that any two complex numbers can be compared through a quotient $r_1/r_2$ of elements of $R$, allowing one small change to act on two positive numbers at the same time. The change is multiplicative $z\mapsto z(1+a)$ but not additive $z\mapsto z+a$ because smallness is relative to the size of $z$.
\begin{definition} \label{def:pair}
A \textit{pair} $(R,\mathfrak{a})$ consists of a subring $R$ of $\mathbb{C}$ whose quotient field is $\mathbb{C}$ and a nonzero ideal $\mathfrak{a}$ of $R$. For a cone $P$ and a pair $(R,\mathfrak{a})$, $P[\mathfrak{a}]=\{0\}\cup\{z\in P^{\times}:\ z(1+\mathfrak{a})\subseteq P\}.$
\end{definition}
\par Thus, $P[\mathfrak{a}]$ is the set of positive numbers whose positivity holds for every change $z\mapsto z(1+a)$ with $a\in\mathfrak{a}$. By Lemma \ref{lem:valring}, $(\mathcal{O}_w,\mathfrak{m}_w)$ and $(\mathcal{O}_w,I^w_\delta)$ are pairs for every nontrivial valuation $w$. For the level valuation, they say beyond the first term and beyond the window of width $\delta$ (Lemma \ref{lem:window}).
\begin{lemma} \label{lem:Pa}
Let $P$ be a cone and $(R,\mathfrak{a})$ a pair.
\begin{enumerate}
\item $P[\mathfrak{a}]$ is a cone, $P\cdot P[\mathfrak{a}]\subseteq P[\mathfrak{a}]$ and
$P[\mathfrak{a}]\,(1+\mathfrak{a})\subseteq P[\mathfrak{a}]$.
\item If $1\in P$, then $P[\mathfrak{a}]=P$ if and only if $1+\mathfrak{a}\subseteq P$.
\item For $p,q\in P[\mathfrak{a}]^{\times}$ there is $d\ne0$ with $p+d$, $p-d$, $q+d$, $q-d\in P$.
\end{enumerate}
\end{lemma}
\begin{proof}
\begin{enumerate}
    \item Let \(x,y\in P[\mathfrak a]\). If \(x+y=0\), then \(x=-y\), and since \(x,y\in P\) we get \(x\in P\cap(-P)=\{0\}\), so \(x=y=0\). For every \(a\in\mathfrak a\),
\[
(x+y)(1+a)=x(1+a)+y(1+a)\in P+P\subseteq P,
\]
so \(x+y\in P[\mathfrak a]\). Thus \(P[\mathfrak a]+P[\mathfrak a]\subseteq P[\mathfrak a]\). Next let \(x,y\in P[\mathfrak a]\) and \(a\in\mathfrak a\). If \(xy=0\), then \(xy\in P[\mathfrak a]\). Otherwise \(x,y\in P^\times\), and \(xy(1+a)\in P\) because \(x(1+a),y(1+a)\in P\) and \(P\cdot P\subseteq P\); also \(xy\ne0\). Hence \(xy\in P[\mathfrak a]\), so \(P[\mathfrak a]\cdot P[\mathfrak a]\subseteq P[\mathfrak a]\). If \(x\in P[\mathfrak a]\cap(-P[\mathfrak a])\), then \(x\in P\) and \(-x\in P\), hence \(x=0\). Thus \(P[\mathfrak a]\cap(-P[\mathfrak a])=\{0\}\), and \(P[\mathfrak a]\) is a cone. Let \(p\in P\) and \(z\in P[\mathfrak a]\). If \(pz=0\), then \(pz=0\in P[\mathfrak a]\). Otherwise \(p,z\in P^\times\), and for every \(a\in\mathfrak a\) we have \(pz(1+a)=p\cdot z(1+a)\in P\cdot P\subseteq P\); hence \(P\cdot P[\mathfrak a]\subseteq P[\mathfrak a]\). Finally, let \(z\in P[\mathfrak a]\) and \(b\in\mathfrak a\). For every \(a\in\mathfrak a\),
\[
z(1+b)(1+a)=z(1+a+b+ab)=z(1+a'')
\]
with \(a''=a+b+ab\in\mathfrak a\), so this lies in \(P\). Thus \(z(1+b)\in P[\mathfrak a]\), and \(P[\mathfrak a](1+\mathfrak a)\subseteq P[\mathfrak a]\).
    \item Assume \(P[\mathfrak a]=P\). Since \(1\in P=P[\mathfrak a]\), $1+\mathfrak a)\subseteq P$. Conversely, assume \(1+\mathfrak a\subseteq P\). We always have \(P[\mathfrak a]\subseteq P\) by definition. Let \(p\in P\). If \(p=0\), then \(p\in P[\mathfrak a]\). If \(p\neq 0\), then \(p\in P^\times\). For every \(a\in\mathfrak a\), we have \(1+a\in P\). Since \(p\in P\) and \(P\cdot P\subseteq P\), $p(1+a)\in P$. Thus \(p(1+\mathfrak a)\subseteq P\), so \(p\in P[\mathfrak a]\). Hence, \( P\subseteq P[\mathfrak a]\), so \(P=P[\mathfrak a]\).
    \item Since \(\mathfrak a\neq\{0\}\), choose \(e\in\mathfrak a^\times\). Because \(R\) has quotient field \(\mathbb C\), the nonzero element \(q/p\in\mathbb C^\times\) can be written as $\frac{q}{p}=\frac{r}{s} $ for some \(r,s\in R\) with \(s\neq 0\). Then, \(r\neq 0\) and \( qs=pr \). Set $a=re\in\mathfrak a$, $b=se\in\mathfrak a$, and $d = pa = qb$. Both \(a\) and \(b\) are nonzero, and \(d\neq 0\). Since \(a,b\in\mathfrak a\), we also have \(-a,-b\in\mathfrak a\). Because \(p,q\in P[\mathfrak a]\), $p+d$, $p-d$, $q+d$, $q-d\in P$.
\end{enumerate}
\end{proof}
\subsection{The Five Conditions}
\par From Section \ref{con}, we know that smallness is measured by a series expansion with a valuation, and positivity depends only on the first term or on a window. In this subsection, each condition is obtained from the previous one by giving up one requirement.
\begin{definition}\label{def:inf}
For a cone $P$, let $\mathcal{I}(P)$ be the set of all $\varepsilon\in\mathbb{C}$ such that $z(1+n\varepsilon)\in P$ for some $z\in P^{\times}$ and all $n\in\mathbb{Z}$.
\end{definition}
\begin{definition}[The five conditions]\label{def:conditions}
For a cone $P\ne\{0\}$ on $\mathbb{C}$:
\begin{enumerate}
\item A: $P=P[\mathfrak{m}_w]$ for some nontrivial valuation $w$ of $\mathbb{C}$.
\item B: $P=P[I^w_\delta]$ for some nontrivial valuation $w$ of $\mathbb{C}$ and some $\delta\ge0$.
\item C: $P=P[\mathfrak{a}]$ for some pair $(R,\mathfrak{a})$.
\item D: $P[\mathfrak{a}]\ne\{0\}$ for some pair $(R,\mathfrak{a})$.
\item E: $\mathcal{I}(P)\ne\{0\}$.
\end{enumerate}
\end{definition}
\par The condition A is for the series expansion whose positivity depends only on the first term. The condition B is for the series expansion whose positivity depends only on a window of fixed width. The condition C indicates that small changes cannot destruct positivity, for some notion of small. The condition D is that they cannot destruct the positivity of at least one number. Condition E says that there are relative infinitesimals; its failure corresponds to the Archimedean case.
\begin{definition}[The six types]\label{def:types}
A directed partial order $P$ on $\mathbb{C}$ is of
\begin{itemize}
\item \textbf{Type 1} if \textup{A} holds (the rule reads only the first term of some series expansion);
\item \textbf{Type 2} if \textup{B} holds but \textup{A} does not (the rule reads a window of fixed width);
\item \textbf{Type 3} if \textup{C} holds but \textup{B} does not (for some pair all positive numbers are insensitive to the small changes, but positivity is not described by a window of one valuation);
\item \textbf{Type 4} if \textup{D} holds but \textup{C} does not (some positive number is insensitive to the small changes, but no single pair makes all positive numbers insensitive);
\item \textbf{Type 5} if \textup{E} holds but \textup{D} does not (there are infinitely small elements, but every positive number is sensitive to arbitrarily small changes);
\item \textbf{Type 6} if \textup{E} does not hold (there are no infinitely small elements).
\end{itemize}
\end{definition}
\par A stratification must achieve the goal that every order is placed in exactly one class. We obtain it because the five conditions are nested. To decide their type, every cone has a condition it satisfies, or satisfies nothing. The proof of their covering is in the rest of this subsection.
\begin{lemma}\label{lem:nested}
For every cone $P\ne\{0\}$ on $\mathbb{C}$,
\[
\textup{A}\Longrightarrow\textup{B}\Longrightarrow\textup{C}\Longrightarrow
\textup{D}\Longrightarrow\textup{E}.
\]
\end{lemma}
\begin{proof}
\par A $\Rightarrow$ B: take $\delta=0$, since $I^{w}_{0} = \mathfrak{m}_{w}$.
\par B $\Rightarrow$ C: $(\mathcal{O}_w,I^w_\delta)$ is a pair by Lemma~\ref{lem:valring}.
\par C $\Rightarrow$ D: $P[\mathfrak{a}]=P\ne\{0\}$.
\par D $\Rightarrow$ E: let $z\in P[\mathfrak{a}]$ and $\varepsilon\in\mathfrak{a}$ be nonzero. For every $n\in\mathbb{Z}$, we have $n\varepsilon\in\mathfrak{a}$, so $z(1+n\varepsilon)\in P$. Hence
$\varepsilon\in\mathcal{I}(P)$.
\end{proof}
\begin{theorem}[Covering Lemma]\label{thm:cover}
For $P \neq 0$, every directed partial order on $\mathbb{C}$ is of exactly one of the Types $1$ to $6$.
\end{theorem}
\begin{proof}
Let $P$ be a directed partial order. Then $P\ne\{0\}$. By Lemma \ref{lem:nested}, if one of the conditions holds, then all later ones hold. Thus, the set of conditions that hold for $P$ is one of $\{A,\dots,E\}$,$\{B,\dots,E\}$, $\{C,D,E\}$, $\{D,E\}$, $\{E\}$, $\emptyset$. These six cases exclude each other and are the types 1, 2, 3, 4, 5, and 6 of Definition \ref{def:types} in this order.
\end{proof}
\subsection{The Six Types}
\par The Covering Lemma holds for the nested conditions, so an order is classified by the first condition it satisfies. The rest of this section explains the types one by one. Types 1 and 2 are described in details; Types 3 and 4 are demonstrated by constructions with examples; Types 5 and 6 are shown by stating what is known because it has not been explored adequately yet.
\subsubsection{Type 1}
\par A basic order $d(S)$ uses the same set $S$ of admissible leading coefficients at every level. Nothing forces it. The general rule reading only the first term possibly uses a different set $L_q$ at each level $q$. We ask which such rules are orders by considering that a sum of two numbers of the same level adds the leading coefficients and a product adds the levels and multiplies the leading coefficients. Next, for the fixed level valuation \(\lv\) of the Puiseux expansion in Section \ref{con}, we describe all cones of depth \(0\) for this fixed valuation; these are a subclass of Type 1. The general condition A allows an arbitrary nontrivial valuation, and the following description does not cover all such valuations.
\par A cone $P$ has \textit{depth} 0 (relative to \(\lv\)) if $P=P[\mathfrak{m}]$ (positivity
depends only on the level and the leading coefficient).
\begin{theorem}\label{thm:type1}
Let $(L_q)_{q\in\mathbb{Q}}$ be a family of subsets of $F^{\times}$ with some empty ones. The rule $P=\{0\}\cup\{z\in\mathbb{C}^{\times}:\ \lc(z)\in L_{\lv(z)}\}$ is a cone if and only if, $\forall q,q'\in\mathbb{Q}$, $L_q+L_q\subseteq L_q$, $L_q\cdot L_{q'}\subseteq L_{q+q'}$. Every cone of depth $0$ arises from exactly one such family. The cone is a directed partial order if and only if the set of levels $\{q:L_q\ne\emptyset\}$ has no lower bound, containing $1$ if and only if $1\in L_0$.
\end{theorem}
\begin{proof}
This is the case $\delta=0$ of Theorems \ref{thm:valid} and \ref{thm:dict} below. For \(\delta=0\), a window is a single coefficient, so the product condition is \(L_qL_{q'}\subseteq L_{q+q'}\), and the sum condition has \(\rho=0\) and becomes \(L_q+L_q\subseteq L_q\); this excludes the case that \(L_q\) contains both \(a\) and \(-a\).
\end{proof}
\noindent Construction: Choose sets $L_q$ satisfying the two conditions of the rule in Theorem \ref{thm:type1}. The four constructions below are chosen by the set of coefficients, the set of levels, a dependence of the coefficients on the level, and a twist by roots of unity.
\begin{enumerate}
\item Let $S$ be a cone on the field $F$ with $1\in S$, and put $L_q=S^{\times}$ $\forall q$. The order is the basic order $d(S)$ of Theorem \ref{thm:basic}. 
\item For $L_q=\mathbb{Q}_{>0}$ with $q\in\mathbb{Z}$ and $L_q=\emptyset$,  only numbers of integer level can be positive.
\item For $L_q=\mathbb{Q}_{>0}$ with $q<0$, $L_0=\mathbb{N}^{\times}$, and $L_q=\emptyset$ with $q>0$, positive numbers have the level at most $0$.
\item $L_q=\chi(q)^{-1}\mathbb{Q}_{>0}$ with $\chi(q)=e^{2\pi iq}$. This is
the order $D^\sigma$ of Example \ref{ex:many}.
\end{enumerate}
\par In each case, the two conditions are checked. We know $1\in L_0$, so the order contains $1$, satisfying Condition A and in Type 1. In particular, all the basic orders of Section \ref{con} are Type 1. To read the number, for example, we take $z=3s^{-1}+is^{-1/2}+\cdots$. Its level is $-1$, and its leading coefficient is $3$. In (1) with $S=\mathbb{Q}_{\ge0}$, in (2) and in (3), the coefficient $3$ is allowed at level $-1$, so $z$ is positive. The later terms are not needed to be looked at.
\subsubsection{Type 2}
\par The next rule reads a fixed number of more terms which is determined by a window. In Type 1, for each level, windows are admissible, and this table is called a card. The question is which cards define orders. There are another two conditions of the direct translations that the product and the sum of positive numbers are positive by windows.
\par A cone $P$ has depth at most $\delta$ if $P=P[I_\delta]$. By Lemma \ref{lem:window}, this means that positivity depends only on the level and the window of width $\delta$.
 
\begin{definition}\label{def:card}
A \textit{card} of width $\delta$ is a family $J=(J_q)_{q\in\mathbb{Q}}$ of subsets of $T^{\times}$, some possibly empty; $J_q$ is the set of \textit{admissible windows at level $q$}, and the \textit{levels} of the card are the $q$ with $J_q\ne\emptyset$. The \textit{rule} of the card is $P_J=\{0\}\cup\{z\in\mathbb{C}^{\times}:\ \win(z)\in J_{\lv(z)}\}.$
\end{definition}
 
\par Every cone of depth at most $\delta$ is the rule of exactly one card: $J_q=\{\win(z):z\in P^{\times},\ \lv(z)=q\}$. Thus, describing these cones is to decide which cards are \textit{valid} (for which cards $P_J$ is a cone).

\par If a rule reads several terms, at the same level, one possibly want that it can use the later terms to accept both a number beginning with $1$ and $-1$. If this can happen, the first terms of two positive numbers can be canceled and their sum can become a completely different level, where the card seems to be useless. The following theorem shows that this never happens. It is the reason why window orders can be described at all and why each of them lies inside a first-term order.
\begin{theorem}\label{thm:sign}
Let $P$ be a cone on $\mathbb{C}$ of depth at most $\delta$. If $z$ and $w$ are positive and have the same level, then $\lc(z)+\lc(w)\ne0$.
\end{theorem}
\begin{proof}
Suppose $\lc(w)=-\lc(z)$. Then $w=-z(1+\eta)$ with $\lv(\eta)>0$. Choose an integer $m\ge1$ with $m\lv(\eta)>\delta$. For $0\le j\le m$, the product $z^{m-j}w^j=(-1)^jz^m(1+\eta)^j$ lies in $P$. Binomial coefficients are natural numbers, so
\[
\sum_{j=1}^m\binom mjz^{m-j}w^j
=z^m\bigl[(1-(1+\eta))^m-1\bigr]
=-z^m\bigl(1-(-\eta)^m\bigr)
\]
lies in $P$. Since $(-\eta)^m\in I_\delta$ and $1+I_\delta$ is a group, depth at most $\delta$ gives $-z^m\in P$. But $z^m\in P^{\times}$, which is a contradiction.
\end{proof}
\begin{example} \label{ex:signnumbers}
Take $\delta=1$, and suppose both $z=1$ and $w=-1+s$ are positive. Then, $2zw+w^2=2(-1+s)+(1-2s+s^2)=-1+s^2$ is positive. A window of width $1$ cannot tell $-1+s^2$ from $-1$, so $-1$ would be positive. However, $1$ is positive, which is a contradiction.
\end{example}
\begin{theorem}\label{thm:valid}
Let $J$ be a card of width $\delta$. The rule $P_J$ is a cone if and only if
\begin{itemize}
\item Products: if $a\in J_q$ and $b\in J_{q'}$, then $ab\in J_{q+q'}$;
\item Sums: if $a\in J_q$ and $b\in J_{q+\rho}$ with $0\le\rho\le\delta$, then $a+s^\rho b\in J_q$.
\end{itemize}
Hence, for the fixed level valuation \(\lv\), the cones on $\mathbb{C}$ of depth at most $\delta$ correspond bijectively to the cards of width $\delta$ satisfying both the conditions of product and sum.
\end{theorem}
\begin{proof} Sufficiency: Assume the condition of the product and the sum. Let $z,z'\in P_J$ be nonzero with levels $q\le q+\rho$ and windows $a,b$. By Lemma \ref{lem:window} and the condition of product, $zz'$ has level $2q+\rho$ and window $ab\in J_{2q+\rho}$. If $\rho>\delta$, then $z+z'$ has level $q$ and window $a$ by Lemma \ref{lem:window}. If $\rho\le\delta$, then $a+s^\rho b\in J_q$ by the condition of the sum. It is a unit, so $z+z'$ has level $q$ and this window. In all cases, $zz'$ and $z+z'$ lie in $P_J$. If $z$ and $-z$ both lay in $P_J$, then $a,-a\in J_q$, and the sum condition with $\rho=0$ would give $0\in J_q$, which is not a unit.
\par Necessity: Let $P_J$ be a cone. It has depth at most $\delta$, and its card is $J$. The condition of the product is derived from Lemma \ref{lem:window}. For the condition of the sum, let $a\in J_q$, $b\in J_{q+\rho}$ with $0\le\rho\le\delta$, and put $z=s^qa$, $z'=s^{q+\rho}b$. The constant term of $a+s^\rho b$ is $a_0$ if $\rho>0$, and $a_0+b_0=\lc(z)+\lc(z')$ if $\rho=0$. In the second case, it is nonzero by Theorem \ref{thm:sign}. By Lemma \ref{lem:window}, the positive number $z+z'$ has level $q$ and window $a+s^\rho b$, which therefore lies in $J_q$. Notice that the sum condition lets the constant of $a+b$ be nonzero when $\rho = 0$, so $J_{q} \cap (-J_{q}) = \emptyset$.
\end{proof}
\par The next theorem is a dictionary showing how to read from the card to determine whether the order is directed, whether $1$ is positive, which constants are positive, and which first-term order contains it.
\begin{theorem}\label{thm:dict}
Let $J$ be a valid card of width $\delta$ with at least one level; $P=P_J$.
\begin{enumerate}
\item The set of levels is closed under addition, and $1\in P$ if and only if $1\in J_0$.
\item $P$ is a directed partial order if and only if the set of levels has no lower bound.
\item Give $L_q=\{a_0:a\in J_q\}$. Then $(L_q)$ satisfies the conditions of multiplication and addition in Type 1, and the first-term order $P^\flat$ of $(L_q)$ contains $P$. It is the smallest cone of depth $0$ containing $P$.
\item $P\cap F=\{0\}\cup\{c\in F^{\times}:c\in J_0\}$.
\end{enumerate}
\end{theorem}
\begin{proof}
\begin{enumerate}
    \item It derives from the condition of Type 2 for the product and $\win(1)=1$.
    \item Suppose the levels have no lower bound. Let $z\ne0$. Choose a positive $b$ with $\lv(b)<\lv(z)-\delta$. Then $z+b=b(1+z/b)$ with $z/b\in I_\delta$, so $z+b$ is positive and $z=(z+b)-b$. Conversely, if all levels are at least $q_0$, then, by Lemma \ref{lem33}, every difference of two elements of $P$ is $0$ or has level at least $q_0$, so $s^{q_0-1}$ is not the difference.
    \item The product condition for \(J\) gives \(L_qL_{q'}\subseteq L_{q+q'}\), since the constant term of \(ab\) is \(a_0b_0\). If \(a,b\in J_q\), then \(a+b\in J_q\) by the sum condition with \(\rho=0\), so \(a_0+b_0\in L_q\), which is the addition condition in Type 1. By Theorem \ref{thm:valid} with $\delta=0$, $P^\flat$ is a cone. It contains $P$, and every cone of depth $0$ containing $P$ contains every number with the level and leading coefficient of an element of $P$. Every order of bounded depth is a first-term order with extra conditions on the later terms of the window.
    \item A nonzero $c\in F$ has level $0$ and window $c$.
\end{enumerate}
\end{proof}
\par We show Type 2 is not empty that some order needs to be seen more than the first term later. The idea is to make the later terms matter in the simplest possible way that a number is positive only if the terms inside its window where their exponents lie in a smaller group $\Lambda$ and their coefficients in a smaller field $K$. The sum and the product still hold the statement, so the card is valid. The condition $(\Lambda,K)\ne(\mathbb{Q},F)$ ensures that the restriction makes sense.
\begin{theorem}[Construction]\label{thm:family}
Let $\Lambda$ be a nonzero subgroup of $\mathbb{Q}$. Let $\delta>0$ be rational. Let $K$ be a subfield of $F$ and $S$ a cone on $K$ with $1\in S$. Assume $(\Lambda,K)\ne(\mathbb{Q},F)$. Let $B$ be the set of truncated series of width $\delta$ whose exponents all lie in $\Lambda$ and whose coefficients all lie in $K$. Give $J_q=\{a\in B:a_0\in S^{\times}\}$ for $q\in\Lambda$ and $J_q=\emptyset$ otherwise. Let $P=P(\Lambda,\delta,K,S)$ be the rule of this card that $z$ is positive if its level lies in $\Lambda$, its leading coefficient lies in $S^{\times}$, and within the window all exponents lie
in $\Lambda$ and all coefficients in $K$. Then $P$ is a directed partial order on $\mathbb{C}$ with $1\in P$ and $P\cap F=S$, which is classified in Type 2.
\end{theorem}
\begin{proof}
\par Firstly, we prove that the card is valid. The ring $B$ is a subring of $T$. The product condition holds because constant terms multiply and $S^{\times} \cdot S^{\times} \subseteq S^{\times}$. For the sum condition, let $a\in J_q$ and $b\in J_{q+\rho}$ with $0\le\rho\le\delta$. Then $\rho\in\Lambda$, so $s^\rho b\in B$, and the constant term of $a+s^\rho b$ is $a_0$ or $a_0+b_0$, which lies in $S^{\times}$. By Theorem \ref{thm:dict}, $P$ is a directed partial order with $1\in P$ because the levels are $\Lambda$ and $1\in J_0$; and $P\cap F=S$.
\par Next, we check how condition B holds. For the level valuation, it is because $P$ has depth at most $\delta$.
\par Finally, we check how condition A fails. If $\Lambda\ne\mathbb{Q}$, choose a rational number $r\notin\Lambda$ with \(0<r\le\delta\); this is possible because if \((0,\delta]\cap\mathbb Q\subseteq\Lambda\), then \(\Lambda\) contains a set of generators of \(\mathbb Q\), hence \(\Lambda=\mathbb Q\). Then \(1+s^r\notin P\). If $\Lambda=\mathbb{Q}$ and $K\ne F$ and  choose $c\in F\setminus K$, then $1+cs^\delta\notin P$. Thus, $1+\mathfrak{m}\not\subseteq P$. Now, suppose $1+\mathfrak{m}_w\subseteq P$ for a nontrivial valuation $w$. If $1+d$ and $1-d$ lie in $P$, then $\lv(d)\ge0$, or they have the same level and leading coefficients $\lc(d)$ and $-\lc(d)$, contradicting Theorem \ref{thm:sign}. Therefore, $\mathfrak{m}_w\subseteq\mathcal{O}$. Let $y\in \mathfrak{m}$ be nonzero. If $y\notin \mathfrak{m}_w$, then $y^{-1}\in\mathcal{O}_w$, and, for a nonzero $a\in \mathfrak{m}_w$, all the numbers $y^{-k}a$ lie in $\mathfrak{m}_w\subseteq\mathcal{O}$. However, $\lv(y^{-k}a)=\lv(a)-k\lv(y)$ is negative for large $k$, so $\mathfrak{m}\subseteq \mathfrak{m}_w$ and $1+ \mathfrak{m}\subseteq P$, which is a contradiction.
\end{proof}
\begin{example}\label{ex:family}
\begin{itemize}
\item $P(\mathbb{Z},1,\mathbb{Q},\mathbb{Q}_{\ge0})$, a two-term order: a positive number has the form $s^n(a_0+a_1s+ \cdots)$ with $n\in\mathbb{Z}$, $a_0\in\mathbb{Q}_{>0}$, $a_1\in\mathbb{Q}$, and no terms strictly inside the window.
\item $P(\tfrac12\mathbb{Z},2,\mathbb{Q}(\sqrt2),S)$ with $S=\{x+y\sqrt2:x,y\in\mathbb{Q}_{\ge0}\}$: There are five coefficients to be read at relative exponents $0,\tfrac12,1,\tfrac32,2$. All lie in
$\mathbb{Q}(\sqrt2)$ and the first in $S^{\times}$.
\end{itemize}
\par In the two-term order, $s^{-2}(3+5s+\cdots)$ is positive, while $s^{-2}(3+s^{1/2}+\cdots)$ is not. It is because of having a term inside the window.
\end{example}
\par For the other expansions, Theorems \ref{thm:type1}--\ref{thm:dict} are stated for the fixed expansion in $s$. A different choice of $T$, of $t$ or of the embedding gives a different expansion where the same theorems apply. We do not determine when two expansions give the same order.
\subsubsection{Type 3}
\par Type 3 consists of orders for which there is a pair \((R,\mathfrak a)\) such that all positive numbers are insensitive to \(\mathfrak a\), but this smallness cannot be measured by a single valuation window. To build one, we use two expansions at once. A function of \(s\) can be expanded near \(s=0\) and also near \(s=1\), and the two expansions see different things. We call a number positive when both expansions are positive. This notion of smallness is given by a pair but not by a single valuation.
\par Set $u=s-1$. Then $F(u)=F(s)$, so every complex number is algebraic over $F(u)$ and has a Puiseux expansion in powers of $u$. Let \(\iota_u\colon\mathbb C\to F\{\{u\}\}\) be the extension to \(\mathbb C\) of \(F(s)=F(u)\to F\{\{u\}\}\) which is the identity on \(F\) and sends \(s\) to \(1+u\); its existence follows as in Lemma \ref{lem32}, since \(\mathbb C\) is algebraic over \(F(u)\) and \(F\{\{u\}\}\) is algebraically closed by Theorem \ref{thm:NP}. We fix this embedding. Let $\lv_u$ and $\lc_u$ be the level and leading coefficient for this second expansion, and $\mathcal{P}_u$ and $\mathfrak{m}_u$ be the corresponding ring and ideal respectively. In the first expansion \(u=-1+s\); in the second \(s=1+u\). Thus \(\lv(s)=1\), \(\lv(u)=0\), \(\lc(u)=-1\), \(\lv_u(u)=1\), \(\lv_u(s)=0\), \(\lc_u(s)=1\).
\begin{theorem}[Construction]\label{thm:type4}
Let $P=\{0\}\cup\{z\in\mathbb{C}^{\times}:\ \lc(z)\in\mathbb{Q}_{>0}\ \text{and}\ \lc_u(z)\in\mathbb{Q}_{>0}\}$. Then $P$ is a directed partial order on $\mathbb{C}$ containing $1$, and it is Type 3.
\end{theorem}
\begin{proof}
\par The set $P$ is the intersection of the two first-term cones with $L_q=\mathbb{Q}_{>0}$ $\forall q$ for each expansion (Theorem \ref{thm:type1}). Thus, it is a cone, and $1\in P$.
\par Let $z\ne0$ and set $b=s^{-N}u^{-2M}$ with integers $N,M\ge1$. In the first expansion $u^{-2M}=(1-s)^{-2M}$ has level $0$ and leading coefficient $1$, so $\lv(b)=-N$ and $\lc(b)=1$. In the second, $s^{-N}=(1+u)^{-N}$ has level $0$ and leading coefficient $1$, so $\lv_u(b)=-2M$ and $\lc_u(b)=1$. For $N>-\lv(z)$ and $2M>-\lv_u(z)$, the number $z+b$ has the same two levels and leading coefficients as $b$. Thus, $b$ and $z+b$ lie in $P$, and $z=(z+b)-b$, suggesting directness.
\par $R=\mathcal{O}\cap\mathcal{O}_u$ and $\mathfrak{a}=\mathfrak{m}\cap\mathfrak{m}_u$. $R$ is a subring and $\mathfrak{a}$ an ideal of $R$. The number $a_0=su$ has $\lv(a_0)=1$ and $\lv_u(a_0)=1$, so $a_0\in\mathfrak{a}$ and $\mathfrak{a}\ne0$. For every $z\in\mathbb{C}$ and $N$ large, $za_0^N\in R$, so $z=(za_0^N)/a_0^N$ lies in the quotient field of $R$, which is $\mathbb{C}$. Thus, $(R,\mathfrak{a})$ is a pair. For $a\in\mathfrak{a}$, the number $1+a$ has level $0$ and leading coefficient $1$ in both expansions, so $1+\mathfrak{a}\subseteq P$, and $P=P[\mathfrak{a}]$ by Lemma \ref{lem:Pa}. Hence, condition C holds.
\par Suppose $P=P[I]$ with $I=I^w_\delta$ for a nontrivial valuation $w$. Then $1+I\subseteq P$. If $d\in I$ had $\lv(d)<0$, then $1+d$ and $1-d$ would have leading coefficients $\lc(d)$ and $-\lc(d)$ without being positive for both of them. Hence, $I\subseteq\mathcal{O}$, and similarly $I\subseteq\mathcal{O}_u$. Consider $x=u/s$, which has $\lv(x)=-1$ and $\lv_u(x)=1$, and choose $a\in I$ nonzero. By Lemma \ref{lem:valring}, $x\in\mathcal{O}_w$ or $x^{-1}\in\mathcal{O}_w$. If $x\in\mathcal{O}_w$, then $x^ka\in I\subseteq\mathcal{O}$ for all $k$, but $\lv(x^ka)=\lv(a)-k<0$ for large $k$. If $x^{-1}\in\mathcal{O}_w$, then $x^{-k}a\in I\subseteq\mathcal{O}_u$, $\forall k$, but $\lv_u(x^{-k}a)=\lv_u(a)-k<0$ for large $k$. Both cases are impossible, so $P$ is of Type 3.
\end{proof}
\par For reading them, the number $\frac{1}{s}$ has $\lc=1$, and $\frac{1}{s}=\frac{1}{1+u}=1-u+\cdots$ has $\lc_u=1$, so it is positive. The number $\frac{1}{s(1-s)}$ has $\lc=1$, but in terms of $u$ it equals $\frac{-1}{u(1+u)}$ with $\lc_u=-1$, so it is not positive, though it is positive for the Type 1 order of the first expansion.
\subsubsection{Type 4}
\par In Types 1--3, every positive number is insensitive to small changes for some pair. Type 4 consists of orders in which this fails for at least one positive number: small changes do not matter all the time. To build one, we stick two opposite parts. For large numbers (of negative level), we keep the first-term rule, which is insensitive. At level \(0\), we put a rigid set where no change is allowed; the natural numbers \(\mathbb N\) are the apparent choice by being closed under addition and multiplication and having gaps between their elements. The two parts can fit together.
\begin{theorem}[Construction]\label{thm:hybrid}
Let $H=\mathbb{N}\ \cup\ \{z\in\mathbb{C}^{\times}:\ \lv(z)<0,\ \lc(z)\in\mathbb{Q}_{>0}\}.$
\begin{enumerate}
\item $H$ is a directed partial order on $\mathbb{C}$ containing $1$.
\item $[-1,1]=\{-1,0,1\}$.
\item $H$ is of Type 4.
\item $H$ is not a lattice order that $0\vee (is^{-1})$ does not exist.
\end{enumerate}
\end{theorem}
A number of negative level is positive when its leading coefficient is a positive rational number. A number of level \(0\) is positive only if it is exactly a natural number. Thus for every \(\delta\), \(2+s^{\delta+1}\) has the same level and the same window of width \(\delta\) as \(2\), but it is not in \(H\). This is why no card of any width describes \(H\).
\begin{proof}
\begin{enumerate}
    \item Firstly, we look at the sum. By Lemma \ref{lem33}, the sum of two positive numbers of negative level has the leading coefficient of the one of lower level, or the sum of the two leading coefficients, which is still a positive rational number. A number of negative level adding a natural number keeps its level and leading coefficient. Secondly, we look at the product. We multiply the leading coefficient. If $z,-z\in H$, then $z=0$. Finally, we check directness. Given $z\ne0$, choose a rational number $q<\min\{\lv(z),0\}$. Then $s^q$ and $z+s^q$ lie in $H$.
    \item If \(\lv(d)<0\), then \(1+d\) and \(1-d\) have leading coefficients \(\lc(d)\) and \(-\lc(d)\), which are not both positive; hence \(d=0\) or \(\lv(d)\ge0\). In the latter case \(1\pm d\) have level at least \(0\), and the elements of \(H\) of level at least \(0\) are exactly the natural numbers. Thus \(1+d=m\) and \(1-d=n\) with \(m,n\in\mathbb N\) and \(m+n=2\), so \(d\in\{-1,0,1\}\). Conversely, if \(d\in\{-1,0,1\}\), then \(1+d,1-d\in\{0,1,2\}\subseteq H\), so \(d\in[-1,1]\). Hence \([-1,1]=\{-1,0,1\}\).
    \item We check how condition D holds and how condition C fails. For condition D, for \(\varepsilon\in \mathfrak{m}\), the number \(s^{-1}(1+\varepsilon)\) has level \(-1\) and leading coefficient \(1\), so \(s^{-1}\in H[\mathfrak{m}]\). For condition C, if \(H=H[\mathfrak{a}]\) for a pair, then \(1+\mathfrak{a}\subseteq H\) by Lemma~\ref{lem:Pa}, so \(\mathfrak{a}\subseteq[-1,1]=\{-1,0,1\}\). However, a nonzero ideal contains all integer multiples of a nonzero element, so it is infinite.
    \item Set \(x=is^{-1}\), which is incomparable with \(0\). Suppose \(h=0\vee x\) exists, so \(h\) and \(h-x\) are positive. If \(\lv(h)>-1\), then \(h-x\) has level \(-1\) and leading coefficient \(-i\), which is not allowed; hence \(\lv(h)\le-1\) and \(\lc(h)\in\mathbb Q_{>0}\). If \(\lv(h)<-1\), choose a rational \(q\) with \(\lv(h)<q<-1\) and put \(h'=h-s^q\). Then \(h'\) and \(h'-x\) have level \(\lv(h)\) and leading coefficient \(\lc(h)\in\mathbb Q_{>0}\), so \(h'\) is an upper bound. But \(h-h'=s^q\in H^\times\), contradicting the minimality of \(h\). If \(\lv(h)=-1\), then \(h-x\) has level \(-1\) and leading coefficient \(\lc(h)-i\notin\mathbb Q\), a contradiction. Therefore \(0\vee(is^{-1})\) does not exist.
\end{enumerate}
\end{proof}
\subsubsection{Type 5 and 6}
\par What remains are the orders in which no positive number is insensitive no matter that the pair is used. They are the types in which a lattice order is not excluded. We divide them once more by determining whether the order is Archimedean. The next theorem shows  condition E and that the relative infinitesimals form a ring.
\begin{theorem}\label{thm:inf}
Let $P$ be a cone on $\mathbb{C}$.
\begin{enumerate}
\item $\mathcal{I}(P)$ is closed under addition, subtraction and multiplication, and $\mathcal{I}(P)\cap\mathbb{Q}=\{0\}$.
\item \(\mathcal{I}(P)=\{0\}\) if and only if \(P\) has no nonzero relative infinitesimal: for all \(z\in P\) and \(y \in \mathbb{C}\), if \(ny\le z\) for all \(n\in\mathbb{Z}\), then \(y=0\).
\end{enumerate}
\end{theorem}
\begin{proof}
\begin{enumerate}
    \item Let $z(1+n\varepsilon)\in P$ and $w(1+n\varepsilon')\in P$, $\forall n\in\mathbb{Z}$ with $z,w\in P^{\times}$. Then $z(1+2n\varepsilon)\cdot w+z\cdot w(1+2n\varepsilon') = 2zw\bigl(1+n(\varepsilon+\varepsilon')\bigr)$, $z(1+n\varepsilon)\cdot w(1+\varepsilon')+z(1-n\varepsilon)\cdot w(1-\varepsilon')=2zw\bigl(1+n\varepsilon\varepsilon'\bigr)$ lie in $P$ for all $n$, and $2zw\in P^{\times}$. Thus, $\varepsilon+\varepsilon'$ and $\varepsilon\varepsilon'$ lie in $\mathcal{I}(P)$. If a nonzero rational number $q$ lie in $\mathcal{I}(P)$, some integer $n$ makes $1+nq$ a negative integer $-m$, and then $-mz\in P$ and $mz\in P$ force $z=0$.
    \item Suppose first that \(\mathcal{I}(P)=\{0\}\). Let \(y,z\in P\) and \(ny\le z\) for all \(n\in\mathbb{Z}\). If \(z=0\), then \(y\le0\) and \(-y\le0\), so \(y\in P\cap(-P)=\{0\}\). If \(z\ne0\), put \(\varepsilon=-y/z\). Then \(z(1+n\varepsilon)=z-ny\in P\) for all \(n\), so \(\varepsilon\in\mathcal{I}(P)=\{0\}\), hence \(y=0\). Conversely, let \(\varepsilon\in\mathcal{I}(P)\), and choose \(z\in P^\times\) with \(z(1+n\varepsilon)\in P\) for all \(n\). Put \(y=-z\varepsilon\). Then \(z(1+n\varepsilon)=z-ny\), so \(ny\le z\) for all \(n\). By the hypothesis, \(y=0\), hence \(\varepsilon=0\). Thus \(\mathcal{I}(P)=\{0\}\).
\end{enumerate}
\end{proof}
\par Thus, Type 6 consists of the directed partial orders with no nonzero relative infinitesimal. Type 5 consists of the orders having relative infinitesimals where, for every pair $(R,\mathfrak{a})$, every positive number loses its positivity under some change $z\mapsto z(1+a)$, $a\in\mathfrak{a}$.
\par \textbf{Does $\mathbb{C}$ have a directed partial order of Type 5 or 6?} We know no construction of an order of either type and cannot prove that none exists. All directed partial orders on $\mathbb{C}$ constructed in this paper are belonged to Types 1 to 4, and we have not determined the types of the orders of \cite{sy11,wyzz26}.
\subsection{Lattice orders}
\par We return to the original question of Birkhoff-Pierce problem. In Theorem \ref{thmnol}, the supremum of $0$ and $i$ failed to exist because an upper bound could be lowered by a number small compared with it with using the expansion. The wiggle Lemma \ref{lemwiggle} isolates what a supremum that a single nonzero change $d$ can be added to and subtracted from both differences $h-x$ and $h-y$. Condition C always provides such $d$ by Lemma \ref{lem:Pa}, so it works for every pair without any expansion.
\begin{theorem}\label{thm:nolattice}
Let $P\ne\{0\}$ be a cone on $\mathbb{C}$ satisfying condition C. Then two complex numbers have a supremum for $P$ only if they are comparable, and $P$ is not a lattice order.
\end{theorem}
\begin{proof}
Let $P=P[\mathfrak{a}]$ for a pair $(R,\mathfrak{a})$. Suppose $x,y$ are incomparable and $h=x\vee y$ exists. By Lemma \ref{lemwiggle}, $p=h-x$ and $q=h-y$ are positive and $[-p,p]\cap[-q,q]=\{0\}$. However, Lemma \ref{lem:Pa} gives a nonzero $d$ in this intersection, so no two incomparable numbers have a supremum. By Lemma \ref{lem21}, for $p\in P^{\times}$, the numbers $0$ and $ip$ are incomparable, so $P$ is not a lattice order. Types 1, 2 and 3 satisfy condition C by Lemma \ref{lem:nested}.
\end{proof}
 
\par When small changes do not matter anymore, any upper bound of two incomparable numbers can be lowered a little and still be an upper bound, so there is no least one. This is the core of Theorem \ref{thmnol}, where the upper bound $h$ of $0$ and $i$ was lowered by $t^{\frac{a}{2}}$. Theorem \ref{thm:nolattice} needs no hypothesis on the set $S$ of the first coefficients. The order $H$ of Type 4 is not a lattice order either by Theorem \ref{thm:hybrid}. We attempt to stratify all the directed partial orders on the complex field, but the question of the existence of lattice order is still remained open.
\section*{Acknowledgements}
\par I would like to express my sincere gratitude to Professor Yuehui Zhang for his guidance and feedback throughout this paper. His insightful comments and suggestions helped me clarify my proof and improve the quality of this work. I am also grateful to him for patiently answering my questions and offering detailed comments on earlier drafts. Any remaining errors are my own.

\end{document}